\documentclass{article} 
\usepackage{amssymb,latexsym}

\newtheorem{thm}{Theorem}[section]

\newtheorem{lem}[thm]{Lemma}

\newtheorem{prop}[thm]{Proposition}
\newtheorem{defn}[thm]{Definition}

\title{\bf\Large Moduli space of Fedosov structures}

\author{ Stanislav Dubrovskiy } 

\date{October 30, 2003}

\begin{document}

\maketitle




\section{Introduction}
Fedosov space is a triple: a manifold $M^{2n}$ with a symplectic
structure $\omega$, and a compatible symmetric connection
$\Gamma$: $(M,\omega,\Gamma)$. Compatibility means that $\omega$
is preserved under geodesic flow of $\Gamma$:
\begin{equation}
\label{eq:compability}
\nabla \omega = 0.
\end{equation}
There is a canonical quantization for these manifolds, see
\cite{GRS} and references therein.

Here we are interested in local invariants of a Fedosov structure.
Namely, we take a space $\mathcal{F}$, of germs of Fedosov
structures at a point, and act on them by local coordinate
changes, that is the group of all diffeomorphisms
$$G:=\mathrm{Diff}(\mathbb{R}^{2n},0)$$ fixing the point. A quotient
of $\mathcal{F}$ by this action is called the moduli space of
Fedosov structures:
$$\mathcal{M}=\mathcal{F}/\mathrm{Diff}(\mathbb{R}^{2n},0).$$
This action can be restricted from space of germs $\mathcal{F}$ to
space of $k$-jets of Fedosov structures, $\mathcal{F}_{k}$, with
corresponding quotient:
$$\mathcal{M}_{k}=\mathcal{F}_{k}/\mathrm{Diff}(\mathbb{R}^{2n},0)
$$
called the moduli space of $k$-jets. We will only work with
generic Fedosov structures. For a generic orbit $\mathcal{O}_k$,
its dimension:
$$\dim\mathcal{\mathcal{O}}_{k}=\textrm{codim }G_{\Phi}
$$
is the codimension of the stabilizer $G_\Phi$ of a generic Fedosov
structure $\Phi$ in $G$. Then we will call
$$\dim\mathcal{M}_{k}=\dim\mathcal{F}_{k}-\dim\mathcal{O}_{k},
$$
and construct the \textit{Poincar\'e series} of $\mathcal{M}$:
$$
p_{\Phi}(t) = \dim \mathcal{M}_{0}+\sum_{k=1}^{\infty}
(\dim\mathcal{M}_{k}- \dim \mathcal{M}_{k-1})t^{k}
$$
\begin{thm}
\label{thm:series} Poncar\'e series coefficients are polynomial in
$k$, and the series has the form:
\begin{eqnarray}
\lefteqn{ \hspace*{-1em}p_{\Phi}(t)  =
\frac{n[8n(2n^2-1)(n+1)+11]}{6}+\frac{n(2n+1)[4n^4+2n^3-6n^2-4n-3]}{3}t+
}
\nonumber\\
& & +(t-t^2)\delta_{2n}^2+{\displaystyle 2n\sum_{k=2}^{\infty}\;\biggl[2{2n+2 \choose
4} {2n+k-1\choose2n-1} - {2n+k+1 \choose 2n-1}\biggr]t^{k}}.
\nonumber
\end{eqnarray}
It represents a rational function.
\end{thm}
{\bf Remark}\ \ This confirms the assertion of Tresse, cf. \cite{T},
that algebras of\\"natural" differential-geometric structures are finitely-generated.\\\\
{\bf Proof}\ \ Postponed until section \ref{sec:Proof}.\\\\
Similar results for other differential-geometric structures were
obtained earlier in \cite{Shmel:struc} and \cite{Dubr}.

To explain significance of Poincar\'e series represented by a rational function, we make the following:\\
{\bf Remark}\ \ 
If a geometric structure is described by a finite number of
functional moduli, then its Poncar\'e series is rational. In
particular, if there are $m$ functional invariants in $n$
variables, then
$$p(t)=\frac{m}{(1-t)^n}$$
Indeed, dimension of moduli spaces of $k$-jets is just the number
of monomials up to the order $k$ in the formal power series of the
$m$ given invariants:$$\dim \mathcal{M}_{k}=m{n+k \choose n}$$ For
more details and slightly more general formulation see {\bf
Theorem 2.1} in \cite{Shmel:mod}.
\section{Action formulas}
As usual, two $C^{\infty}$-functions on $\mathbb{R}^{2n}$ have the
same $k$-jet at a point if their first $k$ derivatives are equal
in any local coordinates.\\

We say that two connections $\nabla$ and $\tilde{\nabla}$ have the
same $k$-jet at 0 if for any two $C^{\infty}$-vector fields $X$,
$Y$ and any $C^{\infty}$-function $f$, the functions
$\nabla_{X}Y(f)$ and $\tilde{\nabla}_{X}Y(f)$ have the same
$k$-jet at $0$. This is equivalent to connection coefficients of
$\nabla$ and $\tilde{\nabla}$ having the same $k$-jet. We denote
by $j^{k}\Gamma$ the $k$-jet of $\Gamma$.\\

There is an action of the group of germs of origin-preserving
diffeomorphisms
$$G=\mathrm{Diff}(\mathbb{R}^{2n},0)$$ on
$\mathcal{F}$ and $\mathcal{F}_{k}$.

For $\varphi \in G$, $(\omega,\nabla)$(or ($\omega,\Gamma)$) $\in
\mathcal{F}$ and $j^{k}\Gamma \in \mathcal{F}_{k}$:
$$\Gamma \mapsto \varphi^{*}\Gamma\,, \quad j^{k}\Gamma \mapsto j^{k}(\varphi^{*}\Gamma)\,,$$
where $${( \varphi^{*}\nabla )}_{X} Y={\varphi}_{*}^{-1}( \nabla_{
\phantom{|}\atop {\varphi}_{*}X  } {\varphi}_{*}Y )$$ Let us
introduce a filtration of $G$ by normal subgroups:
$$ G=G_{1}\rhd G_{2}\rhd G_{3}\rhd \ldots,$$
where $$G_{k}=\{ \,\varphi \in G\ |\
\varphi(x)=x+(\varphi_{1}(x),\ldots\varphi_{n}(x)),\
\varphi_{i}=O(|x|^{k}), \,i=1,\ldots,2n\, \}
$$
The subgroup $G_{k}$ acts trivially on $\mathcal{F}_{p}$ for
$k\geq p+3$.\\
It means that the action of $G$ coincides with that of $G/G_{p+3}$
on each $\mathcal{F}_{p}$.

Now $G/G_{p+3}$ is a finite-dimensional Lie group, which we will
call $K_{p}$ .\\
Denote by $\mathrm{Vect}_{0}(\mathbb{R}^{2n})$ the Lie algebra of
$C^{\infty}$-vector fields, vanishing at the origin. It acts on
$\mathcal{F}$ as follows:
\begin{defn}
For $V\in\mathrm{Vect}_{0}(\mathbb{R}^{2n})$ generating a local
1-parameter subgroup $g^{t}$ of $
\mathrm{Diff}(\mathbb{R}^{n},0)$, \emph{the Lie derivative} of a
connection $\nabla$ in the direction $V$ is a (1,2)-tensor:
$$
\mathcal{L}_{V}\nabla
=\left. \frac{d}{dt}\right|_{t=0}
{g^{t}}^{*}\nabla
$$
\end{defn}
\begin{lem}
\begin{equation}
\label{eq:LD}
(\mathcal{L}_{V}\nabla)(X,Y)=[V,\nabla_{X}Y]-\nabla_{[V,X]}Y-\nabla_{X}[V,Y]
\end{equation}
\end{lem}
{\bf Proof}\quad Below the composition $\circ$ is understood as
that of differential operators acting on functions.
$$
(\mathcal{L}_{V}\nabla)(X,Y)=\left. \frac{d}{dt}\right|_{t=0}
g_{*}^{-t}[ \nabla_{  \phantom{|}\atop { g^{t}_{*}X } } g^{t}_{*}Y
]= \left. \frac{d}{dt}\right|_{t=0} \left[\,(g^{t})^{*}\circ
[\nabla_{  \phantom{|}\atop { g^{t}_{*}X } } g^{t}_{*}Y ] \circ
(g^{-t})^{*}\right]=
$$
$$
\left. \frac{d}{dt}\right|_{t=0} (g^{t})^{*} \circ \nabla_{X}Y +
\nabla_{X}Y \circ \left. \frac{d}{dt}\right|_{t=0} (g^{-t})^{*} +
\nabla_{ \left. \frac{d}{dt}\right|_{t=0}g^{t}_{*}X }Y+
\nabla_{X}\left. \frac{d}{dt}\right|_{t=0}g^{t}_{*}Y =
$$
$$
V\circ \nabla_{X}Y - \nabla_{X}Y \circ V - \nabla_{\left.
\frac{d}{dt}\right|_{t=0}g_{*}^{-t}X}Y
 - \nabla_{X}\left. \frac{d}{dt}\right|_{t=0}g_{*}^{-t}Y =
$$
$$
\mathcal{L}_{V}(\nabla_{X}Y)-\nabla_{\mathcal{L}_{V}X}Y -
\nabla_{X}(\mathcal{L}_{V}Y)
$$
    \hfill$\Box$\\
This defines the action on the germs of connections. Now we can
define the action of $\mathrm{Vect}_{0}(\mathbb{R}^{2n})$ on jets
$\mathcal{F}_{k}$. For $V\in \mathrm{Vect}_{0}(\mathbb{R}^{2n})$:
$$\mathcal{L}_{V}(j^{k}\Gamma)=j^{k}(\mathcal{L}_{V}\Gamma)\ ,$$
where $\Gamma$ on the right
is an arbitrary representative of the $j^{k}\Gamma$ on the left.\\
This is well-defined, since in the coordinate version of
(\ref{eq:LD}):
\begin{equation}
\label{eq:LDcoor} (\mathcal{L}_{V}\Gamma)_{ij}^{l}=
V^{k}\frac{\partial \Gamma_{ij}^{l}}{\partial x^{k}}-
\Gamma_{ij}^{k}\frac{\partial V^{l}}{\partial x^{k}}+
\Gamma_{kj}^{l}\frac{\partial V^{k}}{\partial x^{i}}+
\Gamma_{ik}^{l}\frac{\partial V^{k}}{\partial x^{j}}+
\frac{\partial^{2}V^{l}}{\partial x^{i} \partial x^{j}}
\end{equation}
elements of $k$-th order and less are only coming from
$j^{k}\Gamma$, because $V(0)=0$.
Einstein summation convention in (\ref{eq:LDcoor}) above and further on is assumed.\\
Consequently, the action is invariantly defined. This can also be
expressed as commutativity of the following diagram:
$$
\begin{array}{ccccccccccccc}
j^{0}\mathcal{F}& \longleftarrow & \ldots & \longleftarrow &
j^{k-1}\mathcal{F} & \stackrel{\pi_{k}}{\longleftarrow} &
j^{k}\mathcal{F} & \longleftarrow & \ldots & \longleftarrow &
\mathcal{F} &  &
\\
\downarrow\lefteqn{\mathcal{L}_{V}}& & & &
\downarrow\lefteqn{\mathcal{L}_{V}} & &
\downarrow\lefteqn{\mathcal{L}_{V}} & &
 & & \downarrow\lefteqn{\mathcal{L}_{V}}& &\!\!\!,
\\
j^{0}\Pi & \longleftarrow & \ldots & \longleftarrow & j^{k-1}\Pi &
\stackrel{\pi_{k}}{\longleftarrow} & j^{k}\Pi  & \longleftarrow &
\ldots & \longleftarrow & \Pi & &
\\
\end{array}
$$
where $\pi_{k}$ is projection from $k$-jets onto $(k-1)$-jets,
$\mathcal{F}$ and $\Pi$ denote spaces of germs of connections and
that of (1,2)-tensors respectively, at $0$.

\section{Stabilizer of a generic k-jet}
\label{sec:stab}
[ The following discussion closely mirrors that of section 3 in
\cite{Dubr}. ]
\\Dimensions of stabilizers of generic
k-jets $G_{\Phi}$ are required to find orbit dimensions for orbits
$\mathcal{O}_k$ of generic $k$-jets. The subalgebra generating
$G_{\Phi}$ consists of those
$V\in\mathrm{Vect}_{0}(\mathbb{R}^{2n})$ that
$$
\mathcal{L}_{V}(j^{k}\Phi)=0.
$$
Since $\Phi=(\omega,\Gamma)$, this entails two conditions:
\begin{equation}\label{eq:stab}
\mathcal{L}_{V}(j^{k}\omega)=0\hspace{8ex}\mathcal{L}_{V}(j^{k}\Gamma)=0.
\end{equation}
In the next two sections devoted to finding the stabilizer $G_\Phi$
we assume that $\omega$ is reduced to canonical symplectic form in
Darboux coordinates. In these coordinates compatibility
(\ref{eq:compability}) is written as:
\begin{equation}
\label{eq:compability:coord}
\omega_{i\alpha}\Gamma^\alpha_{kj}=\omega_{j\beta}\Gamma^\beta_{ki},.
\end{equation}
where $\omega=J=\left[\begin{array}{cc}
     \ \ 0       &\textrm{I}\\
-\textrm{I}&    0
 \end{array}\right]
 $, a standard symplectic matrix, cf. \cite{GRS}, p.110.\\
We can introduce grading in homogeneous components on $V$:
$$V=V_{1}+V_{2}+\ldots$$
\begin{center} ( $V_{0}=0$, so that $V$ preserve the origin ) ,
\end{center}
on $\Gamma$:
$$\Gamma=\Gamma_{0}+\Gamma_{1}+\ldots\ ,$$
and on $\omega$:
$$\omega=\omega_0\ ,
$$
where $\omega_0=J$ is a standard symplectic form.\\
Then (\ref{eq:stab}) is rewritten as follows:
$$
\mathcal{L}_{V}(j^{k}\omega)=\mathcal{L}_{V_{1}+V_{2}+\ldots}(\omega_0)=0
$$
$$
\mathcal{L}_{V}(j^{k}\Gamma)=j^{k}\mathcal{L}_{V}(\Gamma)=
j^{k}\mathcal{L}_{V_{1}+V_{2}+\ldots}(\Gamma_{0}+\Gamma_{1}+\ldots+\Gamma_{k}+\ldots)=
$$
$
=\underbrace{\mathcal{L}_{V_{1}}\Gamma_{0}+
\frac{\partial^{2}V_{2}}{\partial x^{2}}}_{\textrm{\footnotesize{0th order}}}+
\underbrace{ \mathcal{L}_{V_{1}}\Gamma_{1}+{\tilde\mathcal{L}}_{V_{2}}\Gamma_{0}+
\frac{\partial^{2}V_{3}}{\partial x^{2}} }_{\textrm{\footnotesize{1st order}}}+
\ldots
$

$\qquad\qquad\qquad\qquad\qquad
\ldots+
\underbrace{\tilde{\mathcal{L}}_{V_{k+1}}\Gamma_{0}+
\tilde{\mathcal{L}}_{V_{k}}\Gamma_{1}+\ldots+\mathcal{L}_{V_{1}}\Gamma_{k}+
\frac{\partial^{2}V_{k+2}}{\partial x^{2}} }_
{\textrm{\footnotesize{k-th order}}}
\quad,
$
$$\mathrm{where}
\quad\left(\frac{\partial^{2}V_{2}}{\partial x^{2}}\right)^{l}_{ij}=
\frac{\partial^{2}V_{2}^{l}}{\partial x^{i}x^{j}}$$
and
$$\quad\tilde{\mathcal{L}}_{V}\Gamma=
\mathcal{L}_{V}\Gamma-\frac{\partial^{2}V}{\partial x^{2}}\ .$$
$\tilde{\mathcal{L}}_{V}\Gamma$ with indexes looks like this:
$$(\tilde{\mathcal{L}}_{V}\Gamma)_{ij}^l=
V^{k}\frac{\partial \Gamma_{ij}^{l}}{\partial x^{k}}-
\Gamma_{ij}^{k}\frac{\partial V^{l}}{\partial x^{k}}+
\Gamma_{kj}^{l}\frac{\partial V^{k}}{\partial x^{i}}+
\Gamma_{ik}^{l}\frac{\partial V^{k}}{\partial x^{j}}\ ,$$
so
$\tilde{\mathcal{L}}_{V}\Gamma$ is just the first 4 terms of
$({\mathcal{L}}_{V}\Gamma)$, cf.(\ref{eq:LDcoor}).\\
The stabilizer condition therefore results in a system:
\begin{equation}\label{sys:SYS}
\left\{
\begin{array}{lr}
\mathcal{L}_{V_{1}}\omega_0=0&\\
\mathcal{L}_{V_{2}}\omega_0=0&\\
\qquad \vdots & \omega-\textrm{part}\\
\mathcal{L}_{V_{k+1}}\omega_0=0&\\
---------------------&\\
\mathcal{L}_{V_{1}}\Gamma_{0}+{\displaystyle\frac{\partial^{2}V_{2}}{\partial x^{2}}} =  0 &\\
\mathcal{L}_{V_{1}}\Gamma_{1}+{\tilde\mathcal{L}}_{V_{2}}\Gamma_{0}+
{\displaystyle\frac{\partial^{2}V_{3}}{\partial x^{2}}}  =  0&\Gamma-\textrm{part}\\
\qquad \vdots  &\\
\mathcal{L}_{V_{1}}\Gamma_{k}+
\tilde{\mathcal{L}}_{V_{2}}\Gamma_{k-1}+\ldots
+\tilde{\mathcal{L}}_{V_{k+1}}\Gamma_{0}+
{\displaystyle\frac{\partial^{2}V_{k+2}}{\partial x^{2}}}  =  0 &
\end{array}
\right.
\end{equation}
Our present goal is finding all $(V_{1},V_{2},\ldots,V_{k+2})$
solving the above system for a generic $\Phi$. Let us start with
the $\Gamma$-part. Assuming $V_{1}$ is arbitrary, we can uniquely
find $V_{2}$ from the first equation, as guaranteed by the
following lemma on mixed derivatives:
\begin{lem}
\label{lem:mixder} Given a family $\{f_{ij}\}_{1\leq i,j\leq N}$
of smooth functions, solution $u$ for the system:
$$
\left\{
\begin{array}{l}
u_{,kl}=f_{kl}\\
1\leq k,l\leq N
\end{array}\right.
$$
(indexes after a comma henceforth will denote differentiations in
corresponding variables) exists if and only if
\begin{equation}
\label{eq:compatibility} \left\{
\begin{array}{l}
f_{ij}=f_{ji}\\
f_{ij,k}=f_{kj,i}
\end{array}\right.
\end{equation}
If $f_{ij}$ are homogeneous polynomials of degree $s\ge0$, then
$u$ can be uniquely chosen as a polynomial of degree $s+2$.
\end{lem}
{\bf Proof}\ \ is a straightforward integration of the right-hand
sides.
\hfill$\Box$\\
Therefore, if we treat highest-order $V_{k}$ in each equation in
$\Gamma$-part of (\ref{sys:SYS}) as an unknown, we see that
various (combinations of) $\mathcal{L}_{V}\Gamma$ must satisfy
(\ref{eq:compatibility}). The first condition is satisfied
automatically since $\Gamma$ is symmetric. The second one gives:
$$(\mathcal{L}_{V}\Gamma)_{ij,p}^{l}=(\mathcal{L}_{V}\Gamma)_{pj,i}^{l}$$
This condition  for the first equation in $\Gamma$-part of
(\ref{sys:SYS}) is satisfied trivially, since $V_{1}$ is of the
first degree in $x$, and $\Gamma_{0}$ is constant. Hence, $V_{2}$
exists and, since it must be of the second degree, is unique.
However, if $k\geq1$ (so there is need for more than one equation)
there is a non-trivial condition on the second equation:
$$\mathcal{L}_{V_1}\Gamma_1 + \widetilde{\mathcal{L}}_{V_2}\Gamma_0 +
\frac{\partial^2 V_3}{\partial x^{2}} = 0$$ It follows from the
Lemma\,\ref{lem:mixder} that for the existence of $V_3$ it is
necessary ( and sufficient ) to have the following condition:
\begin{equation}
\label{eq:3.7} (\mathcal{L}_{V_1}\Gamma_1 +
\widetilde{\mathcal{L}}_{V_2}\Gamma_0)_{ij,p}=
(\mathcal{L}_{V_1}\Gamma_1 +
\widetilde{\mathcal{L}}_{V_2}\Gamma_0)_{pj,i}\ ,\ i<p
\end{equation}
Outside exceptional dimension two this condition fails for a
generic connection unless $V_1=0$. In other words (\ref{eq:3.7}),
considered as a condition on $V_1$ implies $V_1=0$ ( and hence
$V_2=V_3=...=0$ ). The rest is a proof of this assertion.

Let us consider (\ref{eq:3.7}) as a linear homogeneous system on
the components of $V_1$. We will present a Fedosov structure for
which (\ref{eq:3.7}) is non-degenerate. Since nondegeneracy is an
open condition on the space $\mathcal{F}_k$, the same (namely
non-degeneracy and resulting trivial solution $V=0$ for the
stabilizer) would hold for a generic structure.

Since symplectic part of the structure is fixed by our choice to
work in symplectic coordinates, we need only to present the
corresponding connection 1-jet. In this 1-jet we choose to have
$\Gamma_0=0$, which simplifies (\ref{eq:3.7}) to:
$$(\mathcal{L}_{V_1}\Gamma_1)_{ij,p}^{l}=(\mathcal{L}_{V_1}\Gamma_1 )_{pj,i}^{l}$$
Let us expand it using (\ref{eq:LDcoor}):
\\[2.5ex]
$(\Gamma_{1 ij,k}^{l}-\Gamma_{1 kj,i}^{l})V_{1 ,p}^{k}+ (\Gamma_{1
kj,p}^{l}-\Gamma_{1 pj,k}^{l})V_{1 ,i}^{k}+$
\begin{equation}
\label{eq:3.9}
\qquad (\Gamma_{1 pj,i}^{k}-\Gamma_{1
ij,p}^{k})V_{1 ,k}^{l}+ (\Gamma_{1 ik,p}^{l}-\Gamma_{1
pk,i}^{l})V_{1 ,j}^{k}=0\ ,\ i<p
\end{equation}
Recall that summation over repeated indexes above is assumed.

In local symplectic coordinates:
\[
V_{1}^{k}=\sum_{s=1}^{n}b_{s}^{k}x^{s} \ ,\qquad \Gamma_{1
ij}^{l}=\sum_{m=1}^{n}c_{ij}^{lm}x^{m}, \quad
c_{ij}^{lm}=c_{ji}^{lm}\textrm{  (connection is symmetric),}
\]
and (\ref{eq:3.9}) becomes the system on $b_{s}^{k}$:
\begin{equation}
\label{eq:3.10}
(c_{ij}^{lk}-c_{kj}^{li})b_{p}^{k}+
(c_{kj}^{lp}-c_{pj}^{lk})b_{i}^{k}+
(c_{pj}^{ki}-c_{ij}^{kp})b_{k}^{l}+
(c_{ik}^{lp}-c_{pk}^{li})b_{j}^{k}=0\ ,\ i<p
\end{equation}
The requirement
(\ref{eq:compability}) on $\Gamma$ to be a symplectic connection is passed through to each of its
homogeneous components $\Gamma_\textbf{k}$ as the following symmetry condition:
$$
\omega_{i\alpha}\Gamma^\alpha_{\textbf{k}jl}=\omega_{l\alpha}\Gamma^\alpha_{\textbf{k}ji}\ ,
$$
that can be thought of as `$\Gamma$ with lowered indexes' is completely symmetric
(cf. the discussion on p.110 (especially equation (1.5)) in \cite{GRS}).\\
Of course $\omega_{i\alpha}$ in our setting is just the standard
symplectic matrix. Another way to think about it in terms of e.g.
coefficients of $\Gamma_1$ is that they `form a symplectic
matrix', namely $\Gamma_1\in \texttt{sp}(2n)$ in the left upper
and lower indexes:
\begin{eqnarray}
\label{eq:2.7}
c^{im}_{Jk} & = & c^{\bar{J}m}_{\bar{i}k}\nonumber\\
c^{Im}_{jk} & = & c^{\bar{j}m}_{\bar{I}k}\\
c^{Im}_{Jk} & = & -c^{\bar{J}m}_{\bar{I}k}\ \Rightarrow\ \ c^{im}_{jk}\ =\ -c^{\bar{j}m}_{\bar{i}k} \nonumber
\end{eqnarray}
$$\forall m, k \in[1,\ldots 2n],\ i,j \in[1\ldots n],\ I,J \in [n+1,\ldots 2n],\ \bar{i}=i+n,\bar{I}=I-n$$
We will also consider only such $\Gamma_1$ that
\begin{equation}
\label{3.11} c_{ij}^{lp}\neq0
\textrm{  only if }\{i,j,l,p\}=\{\alpha,\beta,\gamma\}\,,\ \alpha\neq\beta,\,\beta\neq\gamma,\,\alpha\neq\gamma
\end{equation}
In other words nonzero coefficients may only occur among those
with indexing set consisting of three distinct numbers, and must
be zero otherwise.\\
%
%
%
%
%

Let us now turn to $\omega$-part of (\ref{sys:SYS}).
The fact that $V_1$ preserves $\omega$ implies that it is
hamiltonian, so its coefficients:
$$
b=J\nabla H\,
$$
for some hamiltonian ${\displaystyle H=\frac{1}{2}\sum_{i,j=1}^{n} h_{ij}x^i x^j}$, homogeneous\\
second degree polynomial,
i.e.:
\begin{equation}
\label{comp:Gamma1}
b^k_s=\pm h_{k\pm n,s}\,
\end{equation}
where the sign that makes sense applies:
\begin{displaymath}
\left[
\begin{array}{rcl}
+ & for & k\leq n\\
- & for & k> n
\end{array}
\right.
\end{displaymath}
Our choice of $\Gamma_0=0$ implies $V_2=0$. Hence second equation in $\omega$-part of (\ref{sys:SYS})
is satisfied trivially.

(\ref{3.11}), coupled with compatibility conditions (\ref{eq:2.7}) leaves us only five types
of possibly non-zero coefficients:
$$
c^{\alpha m}_{\alpha k},\ m\neq k,\ \{m,k\}\cap\{\alpha,\bar{\alpha}\}=\emptyset\ ,
\qquad
c^{\alpha k}_{\beta k},\ k\not\in \{\alpha,\beta,\bar{\alpha},\bar{\beta}\},\ \alpha\neq \beta\ ,
$$
\nopagebreak[4]
$$
c^{\alpha \alpha}_{\bar{\alpha} k},\ c^{\alpha \bar{\alpha}}_{\bar{\alpha} k},\
c^{\alpha k}_{\bar{\alpha} \bar{\alpha}},\ k\not\in \{\alpha,\bar{\alpha}\}\ .
$$\\
We will now specify (\ref{eq:3.10}) to different particular $\{i,j,l,p\}$ .\\
i) $i=\bar{P}$, $j=l\neq i$
$$
(c_{ij}^{jk}-c_{kj}^{ji})b_{\bar{i}}^{k}+
(c_{kj}^{j\bar{i}}-c_{\bar{i}j}^{jk})b_{i}^{k}+
(c_{\bar{i}j}^{ki}-c_{ij}^{k\bar{i}})b_{k}^{j}+
(c_{ik}^{j\bar{i}}-c_{\bar{i}k}^{ji})b_{j}^{k}=0
$$
There are 3 distinct indexes present in each coefficient in the above equation.
If it seems that there are only 2, we must use their symmetries (\ref{eq:2.7}), to explicitly present all three.
For example, the first coefficient $c_{ij}^{jk}=-c_{\bar{j}j}^{\bar{i}k}$, and similar for other coefficients.
That implies that in each summation the dummy index $k$ has to turn into one of the fixed ones, e.g.
into $j$, $\bar{i}$ or $\bar{j}$ in the first term:
$$(c_{ij}^{jj}-c_{jj}^{ji})b_{\bar{i}}^{j}+
(c_{ij}^{j\bar{i}}-c_{\bar{i}j}^{ji})b_{\bar{i}}^{\bar{i}}+
(c_{ij}^{j\bar{j}}-c_{\bar{j}j}^{ji})b_{\bar{i}}^{\bar{j}}=
(c_{ij}^{jj}-c_{jj}^{ji})h_{\bar{j}\bar{i}}-
(c_{ij}^{j\bar{i}}-c_{\bar{i}j}^{ji})h_{i\bar{i}}-
(c_{ij}^{j\bar{j}}-c_{\bar{j}j}^{ji})h_{j\bar{i}}
$$
Note that even though indexes do repeat in the first term above, the summation convention does not apply,
because the designated summation dummy $k$ is absent!
After the similar work is done for the remaining 3 terms, and many cancellations
(due to symmetries (\ref{eq:2.7})),the original equation simplifies to:
$$
c_{\bar{i}\bar{i}}^{ji}h_{ij}+c_{\bar{i}j}^{ii}h_{i\bar{j}}+
c^{j\bar{i}}_{ii}h_{\bar{i}j}-c_{ij}^{\bar{i}\bar{i}}h_{\bar{j}\bar{i}}=0\qquad i\neq j\ .
$$
In the similar manner we obtain the next four equations:\\
ii) $i=\bar{P}$, $j=\bar{L}\neq i$
$$
2(c_{ji}^{j\bar{i}}-c_{j\bar{j}}^{ii})h_{jj}+(c_{j\bar{i}}^{\bar{j}\bar{j}}-2c_{\bar{i}j}^{ii})h_{ij}+
c_{jj}^{\bar{j}\bar{i}}h_{i\bar{j}}+(2c_{ii}^{\bar{j}\bar{i}}-c_{ij}^{\bar{j}\bar{j}})h_{\bar{i}j}-
c_{jj}^{\bar{j}i}h_{\bar{i}\bar{j}}=0\qquad i\neq j\ .
$$\\
iii) $\bar{I}=j(=:i)$, $l=\bar{P}$, $I<P\Rightarrow j(=i)<l$
$$(c_{\bar{i}\bar{i}}^{l\bar{l}}-c_{\bar{l}\bar{i}}^{l\bar{i}})h_{ii}-
c_{\bar{l}\bar{l}}^{l\bar{i}}h_{il}
+(2c_{ii}^{l\bar{i}}-c_{\bar{l}\bar{l}}^{\bar{i}l})h_{\bar{i}\bar{l}}+
(c_{\bar{l}i}^{li}-c_{ii}^{l\bar{l}})h_{\bar{i}\bar{i}}+(2c_{li}^{l\bar{i}}-c_{\bar{i}i}^{ll}-c_{il}^{i\bar{l}})
h_{\bar{l}\bar{l}}=0\qquad i<l\ .
$$
\quad $\bar{I}=l$, $j=\bar{P}$, $I<P\Rightarrow l<j$
$$(c_{\bar{j}\bar{j}}^{l\bar{l}}-c_{\bar{l}\bar{j}}^{l\bar{j}})h_{jj}-
c_{\bar{l}\bar{l}}^{l\bar{j}}h_{jl}
+(2c_{jj}^{l\bar{j}}-c_{\bar{l}\bar{l}}^{\bar{j}l})h_{\bar{j}\bar{l}}+
(c_{\bar{l}j}^{lj}-c_{jj}^{l\bar{l}})h_{\bar{j}\bar{j}}+(2c_{lj}^{l\bar{j}}-c_{\bar{j}j}^{ll}-c_{jl}^{j\bar{l}})
h_{\bar{l}\bar{l}}=0\qquad j>l\ .
$$
Since these two equations are the same modulo changing $i$ into $j$,
we can keep just the last equation, but for $j\neq l$.
And finally we rewrite it in $i$ and $j$ in conformity with others:
$$
(c_{\bar{i}\bar{i}}^{j\bar{j}}-c_{\bar{j}\bar{i}}^{j\bar{i}})h_{ii}-
c_{\bar{j}\bar{j}}^{j\bar{i}}h_{ij}
+(2c_{ii}^{j\bar{i}}-c_{\bar{j}\bar{j}}^{\bar{i}j})h_{\bar{i}\bar{j}}+
(c_{\bar{j}i}^{ji}-c_{ii}^{j\bar{j}})h_{\bar{i}\bar{i}}+(2c_{ji}^{j\bar{i}}-c_{\bar{i}i}^{jj}-c_{ij}^{i\bar{j}})
h_{\bar{j}\bar{j}}=0\qquad i\neq j\ .
$$
iv) $i=l$, $J(=:\bar{j}\ )=P\neq \bar{i}$
$$(c_{\bar{j}\bar{j}}^{i\bar{i}}-c_{\bar{i}\bar{j}}^{i\bar{j}})h_{ii}
+(c_{\bar{j}\bar{j}}^{\bar{i}i}-c_{i\bar{j}}^{\bar{i}\bar{j}})h_{\bar{i}\bar{i}}
+(c_{i\bar{j}}^{ij}+c_{ij}^{i\bar{j}}-2c_{j\bar{j}}^{ii})h_{\bar{j}\bar{j}}
+(2c_{\bar{i}\bar{j}}^{ii}-c_{\bar{j}\bar{j}}^{ij})h_{i\bar{j}}
+c_{\bar{j}\bar{j}}^{ji}h_{\bar{i}j}=0\qquad i\neq j\ .
$$\\
v) $i=\bar{L}$, $J(=:\bar{j}\ )=P\neq \bar{i}$
$$-c_{\bar{j}\bar{j}}^{ji}h_{ij}
+2(c_{i\bar{j}}^{\bar{i}\bar{j}}-c_{\bar{j}\bar{j}}^{\bar{i}i})h_{i\bar{i}}
+(c_{\bar{j}\bar{j}}^{\bar{i}i}-c_{i\bar{j}}^{\bar{i}\bar{j}})h_{j\bar{j}}
+c_{ij}^{\bar{j}\bar{j}}h_{\bar{j}\bar{j}}
+(c_{i\bar{j}}^{jj}-c_{i\bar{j}}^{\bar{i}\bar{i}}+c_{\bar{j}i}^{\bar{j}i}-c_{i\bar{i}}^{\bar{j}\bar{j}})
h_{i\bar{j}}=0\qquad i\neq j\ .
$$
In each of the equations i)-v) above we are free to interchange $i$ with $j$ to obtain another five:
i'), ii'), iii'), iv') and v'). Thus we obtain the system of ten equations for ten variables:
$h_{ij}$, $h_{i\bar{j}}$, $h_{\bar{i}j}$, $h_{\bar{i}\bar{j}}$, $h_{ii}$, $h_{\bar{i}\bar{i}}$, $h_{jj}$,
$h_{\bar{j}\bar{j}}$, $h_{i\bar{i}}$, and $h_{j\bar{j}}$. However, the last two variables are only found in
the equation v), which allows us to consider first four equations and their `primes' i), i'),... iv')
as an 8 X 8 system for the first eight variables:\\\\
$\hspace*{-4.8em}
h_{ij}\hspace*{4.5em} h_{i\bar{j}}\hspace*{4.5em} h_{\bar{i}j}\hspace*{4.5em} h_{\bar{i}\bar{j}} \hspace*{4.7em}
h_{ii} \hspace*{5.1em} h_{\bar{i}\bar{i}} \hspace*{6em} h_{jj} \hspace*{5.8em} h_{\bar{j}\bar{j}}
$
$$
\hspace*{-8.5em}
\left(
\begin{array}{cccccccc}
-c_{\bar{i}\bar{i}}^{ji} & -c_{\bar{i}j}^{ii} & -c^{j\bar{i}}_{ii} & c_{ij}^{\bar{i}\bar{i}} \\
 & & & & & & & \\
-c_{\bar{j}\bar{j}}^{ij} & -c^{i\bar{j}}_{jj} & -c_{\bar{j}i}^{jj} & c_{ji}^{\bar{j}\bar{j}} \\
 & & & & & & & \\
(c_{\bar{i}j}^{\bar{j}\bar{j}}-2c_{\bar{i}j}^{ii}) & c_{jj}^{\bar{j\bar{i}}} &
(2c_{ii}^{\bar{j}\bar{i}}-c_{ij}^{\bar{j}\bar{j}}) & -c_{jj}^{\bar{j}i} &
& & 2(c_{ji}^{j\bar{i}}-c_{j\bar{j}}^{ii}) & \\
 & & & & & & & \\
(c_{i\bar{j}}^{\bar{i}\bar{i}}-2c_{\bar{j}i}^{jj}) & (2c_{jj}^{\bar{i}\bar{j}}-c_{ji}^{\bar{i}\bar{i}}) &
c_{ii}^{\bar{i\bar{j}}} & -c_{ii}^{\bar{i}j} &
2(c_{ij}^{i\bar{j}}-c_{i\bar{i}}^{jj}) & &  & \\
 & & & & & & & \\
-c_{\bar{j}\bar{j}}^{j\bar{i}} & & & (2c_{ii}^{j\bar{i}}-c_{\bar{j}\bar{j}}^{\bar{i}j}) &
(c_{\bar{i}\bar{i}}^{j\bar{j}}-c_{\bar{j}\bar{i}}^{j\bar{i}}) &
(c_{\bar{j}i}^{ji}-c_{ii}^{j\bar{j}}) & & (2c_{ji}^{j\bar{i}}-c_{\bar{i}i}^{jj}-c_{ij}^{i\bar{j}}) \\
 & & & & & & & \\
-c_{\bar{i}\bar{i}}^{i\bar{j}} & & & (2c_{jj}^{i\bar{j}}-c_{\bar{i}\bar{i}}^{\bar{j}i}) &
 & (2c_{ij}^{i\bar{j}}-c_{\bar{j}j}^{ii}-c_{ji}^{j\bar{i}})
 & (c_{\bar{j}\bar{j}}^{i\bar{i}}-c_{\bar{i}\bar{j}}^{i\bar{j}}) &(c_{\bar{i}j}^{ij}-c_{jj}^{i\bar{i}}) \\
 & & & & & & & \\
 & (2c_{\bar{i}\bar{j}}^{ii}-c_{\bar{j}\bar{j}}^{ij}) & c_{\bar{j}\bar{j}}^{ji} & &
 (c_{\bar{j}\bar{j}}^{i\bar{i}}-c_{\bar{i}\bar{j}}^{i\bar{j}}) &
 (c_{\bar{j}\bar{j}}^{\bar{i}i}-c_{i\bar{j}}^{\bar{i}\bar{j}}) & &
 (c_{i\bar{j}}^{ij}+c_{ij}^{i\bar{j}}-2c_{j\bar{j}}^{ii}) \\
 & & & & & & & \\
 & c_{\bar{i}\bar{i}}^{ij} & (2c_{\bar{j}\bar{i}}^{jj}-c_{\bar{i}\bar{i}}^{ji}) & &
  &(c_{j\bar{i}}^{ji}+c_{ji}^{j\bar{i}}-2c_{i\bar{i}}^{jj})
  & (c_{\bar{i}\bar{i}}^{j\bar{j}}-c_{\bar{j}\bar{i}}^{j\bar{i}})
 & (c_{\bar{i}\bar{i}}^{\bar{j}j}-c_{j\bar{i}}^{\bar{j}\bar{i}}) \\
\end{array} \right)$$\\
We set:
$$
2c_{ii}^{j\bar{i}}-c_{\bar{j}\bar{j}}^{\bar{i}j}
(=2c_{i\bar{j}}^{\bar{i}\bar{i}}-c_{\bar{j}i}^{jj})=0\ ,
\quad c_{\bar{i}\bar{j}}^{ii}(=c_{\bar{i}\bar{i}}^{ji})=0\ ,\quad
c_{**}^{*\alpha}=0\ ,*\in\{i,\bar{i}\}\ ,\ \alpha\not\in \{i,\bar{i}\}\ .
$$
This doesn't completely separate the system, but it does annihilate the lower left block.
Consider the lower right block:\\\\
$\hspace*{6em}
h_{ii} \hspace*{5.4em} h_{\bar{i}\bar{i}} \hspace*{5.4em} h_{jj} \hspace*{5.4em} h_{\bar{j}\bar{j}}
$
$$
\left(
\begin{array}{cccc}
(c_{\bar{i}\bar{j}}^{i\bar{j}}-c_{\bar{j}\bar{i}}^{j\bar{i}}) &
(c_{\bar{j}i}^{ji}+c_{i\bar{j}}^{\bar{i}\bar{j}}) & 0 & (2c_{ji}^{j\bar{i}}-c_{i\bar{j}}^{ij}-c_{ij}^{i\bar{j}}) \\
 & & & \\
 0 & (2c_{ij}^{i\bar{j}}-c_{j\bar{i}}^{ji}-c_{ji}^{j\bar{i}})
 & (c_{\bar{j}\bar{i}}^{j\bar{i}}-c_{\bar{i}\bar{j}}^{i\bar{j}}) & (c_{\bar{i}j}^{ij}+c_{j\bar{i}}^{\bar{j}\bar{i}}) \\
 & & & \\
 (c_{\bar{j}\bar{i}}^{j\bar{i}}-c_{\bar{i}\bar{j}}^{i\bar{j}}) &
 -(c_{\bar{j}i}^{ji}+c_{i\bar{j}}^{\bar{i}\bar{j}}) & 0 &
 (2c_{\bar{j}\bar{i}}^{\bar{j}i}+c_{ij}^{i\bar{j}}+c_{i\bar{j}}^{ij}) \\
 & & & \\
 0 & (2c_{\bar{i}\bar{j}}^{\bar{i}j}+c_{j\bar{i}}^{ji}+c_{ji}^{j\bar{i}})
  & (c_{\bar{i}\bar{j}}^{i\bar{j}}-c_{\bar{j}\bar{i}}^{j\bar{i}})
 & -(c_{\bar{i}j}^{ij}+c_{j\bar{i}}^{\bar{j}\bar{i}}) \\
\end{array} \right)$$\\
It has ten independent coefficients, which we can set as follows:
$$
c_{\bar{i}\bar{j}}^{i\bar{j}}=
\left[
\begin{array}{rl}
+1 & i>j\\
-1 & i<j
\end{array}
\right.
\ ,\quad
c_{ji}^{j\bar{i}}=c_{ij}^{i\bar{j}}=2\ ,\textrm{ and the rest all equal to 1. }
$$
It is clearly non-degenerate, hence we just need to show that the upper left block also can be chosen non-degenerate:\\\\
$\hspace*{-0.64em}
h_{ij}\hspace*{3.7em} h_{i\bar{j}}\hspace*{4.5em} h_{\bar{i}j}\hspace*{2.4em} h_{\bar{i}\bar{j}}
$
$$\hspace*{-2.8em}
\left(
\begin{array}{cccc}
0 & -c_{\bar{i}j}^{ii} & -c^{j\bar{i}}_{ii} & c_{ij}^{\bar{i}\bar{i}} \\
 & & & \\
0 & -c^{i\bar{j}}_{jj} & -c_{\bar{j}i}^{jj} & c_{ji}^{\bar{j}\bar{j}} \\
 & & & \\
-3c_{j\bar{i}}^{\bar{j}\bar{j}} & 0 &
(2c_{ii}^{\bar{j}\bar{i}}-c_{ij}^{\bar{j}\bar{j}}) & 0 \\
 & & & \\
-3c_{i\bar{j}}^{\bar{i}\bar{i}} & (2c_{jj}^{\bar{i}\bar{j}}-c_{ji}^{\bar{i}\bar{i}}) &
0 & 0 \\
\end{array} \right)=
\left(
\begin{array}{cccc}
0 & -c_{\bar{i}j}^{ii} & \frac{1}{2}c^{\bar{j}i}_{jj} & c_{ij}^{\bar{i}\bar{i}} \\
 & & & \\
0 & \frac{1}{2}c^{i\bar{j}}_{jj} & -c_{\bar{j}i}^{jj} & c_{ji}^{\bar{j}\bar{j}} \\
 & & & \\
-\frac{3}{2}c_{\bar{i}j}^{ii} & 0 &
(2c_{ij}^{\bar{i}\bar{i}}-c_{ji}^{\bar{j}\bar{j}}) & 0 \\
 & & & \\
-\frac{3}{2}c_{\bar{j}i}^{jj} & (2c_{ji}^{\bar{j}\bar{j}}-c_{ij}^{\bar{i}\bar{i}}) &
0 & 0 \\
\end{array} \right)$$
There are four independent coefficients in this system:
$$
c_{\bar{i}j}^{ii}\,,\ c^{\bar{j}i}_{jj}\,,\ c_{ij}^{\bar{i}\bar{i}}\ \textrm{ and }\ c_{ji}^{\bar{j}\bar{j}}\,.
$$
Setting them all equal to 1 achieves non-degeneracy for this block, and for the 8 X 8 system.
This leaves us equations v) and v'), which reduce to this 2 X 2 system for the remaining two
unknowns $h_{i\bar{i}}$ and $h_{j\bar{j}}$:
$$
\left(
\begin{array}{cc}
2(c_{i\bar{j}}^{\bar{i}\bar{j}}+c^{ji}_{\bar{j}i}) & -(c_{\bar{j}i}^{ji}+c_{i\bar{j}}^{\bar{i}\bar{j}}) \\
                                                   &                                                   \\
-(c^{ij}_{\bar{i}j}+c_{j\bar{i}}^{\bar{j}\bar{i}}) & 2(c_{j\bar{i}}^{\bar{j}\bar{i}}+c_{\bar{i}j}^{ij})\\
\end{array} \right)
=
\left(
\begin{array}{cc}
 4 & -2 \\
   &    \\
-2 & 4  \\
\end{array} \right)$$
Since we required 3 distinct indexes for our non-zero coefficients,
this method is only good for dimensions 4 or larger.
\begin{prop}\label{stab:highdim}
The stabilizer of a k-jet of a generic connection for $n\geq 4$ is:
$G_{1}/G_{2}$
 for $k=0$
, and $0$
 for $k\geq 1$ .
\end{prop}
The lowest dimension 2 has to be treated separately.
\section{Exceptional dimension 2}
In this case the stabilizer of the 1-jet is non-trivial (it has dimension one),
stabilizers of the higher jets are all trivial.\\
Since the general method of previous section fails here, we must reconsider (\ref{eq:3.10})
with $i=1$ and $p=2$:
$$
(c_{1j}^{lk}-c_{kj}^{l1})b_{2}^{k}+
(c_{kj}^{l2}-c_{2j}^{lk})b_{1}^{k}+
(c_{2j}^{k1}-c_{1j}^{k2})b_{k}^{l}+
(c_{1k}^{l2}-c_{2k}^{l1})b_{j}^{k}=0
$$
Summing over two indexes, we obtain:
$$
(c_{2j}^{11}-c_{1j}^{12})b_{1}^{l}+(c_{2j}^{21}-c_{1j}^{22})b_{2}^{l}+
(c_{11}^{l2}-c_{21}^{l1})h_{2j}-(c_{12}^{l2}-c_{22}^{l1})h_{1j}=0
$$
Varying pair $(ij)$ we obtain next 4 equations on 3 variables $h_{11}$, $h_{12}$ and $h_{22}$:
$$
(c_{21}^{11}-c_{11}^{12})h_{21}+(c_{21}^{21}-c_{11}^{22})h_{22}+
(c_{11}^{12}-c_{21}^{11})h_{21}-(c_{12}^{12}-c_{22}^{11})h_{11}=0
$$
$$
\hspace*{-1.68ex}-(c_{21}^{11}-c_{11}^{12})h_{11}-(c_{21}^{21}-c_{11}^{22})h_{12}+
(c_{11}^{22}-c_{21}^{21})h_{21}-(c_{12}^{22}-c_{22}^{21})h_{11}=0
$$
$$
(c_{22}^{11}-c_{12}^{12})h_{21}+(c_{22}^{21}-c_{12}^{22})h_{22}+
(c_{11}^{12}-c_{21}^{11})h_{22}-(c_{12}^{12}-c_{22}^{11})h_{12}=0
$$
$$
\hspace*{-1.68ex}-(c_{22}^{11}-c_{12}^{12})h_{11}-(c_{22}^{21}-c_{12}^{22})h_{12}+
(c_{11}^{22}-c_{21}^{21})h_{22}-(c_{12}^{22}-c_{22}^{21})h_{12}=0
$$
The coefficient matrix of the system is this:
\\
$\hspace*{8.5em}
h_{11}\hspace*{7.5em} h_{12}\hspace*{7em} h_{22}
$
$$
\left(
\begin{array}{ccc}
-(c_{12}^{12}-c_{22}^{11}) & & (c_{21}^{21}-c_{11}^{22})\\
 & & \\
-(c_{21}^{11}-c_{11}^{12}+c_{12}^{22}-c_{22}^{21}) & -2(c_{21}^{21}-c_{11}^{22}) & \\
 & & \\
 & 2(c_{22}^{11}-c_{12}^{12}) &
(c_{22}^{21}-c_{12}^{22}+c_{11}^{12}-c_{21}^{11})\\
 & & \\
-(c_{22}^{11}-c_{12}^{12}) & & (c_{11}^{22}-c_{21}^{21})\\
\end{array} \right)
$$
Setting
$$
a:=c_{11}^{12}-c_{21}^{11}\ ,\quad b:=c_{12}^{12}-c_{22}^{11}\ ,\quad \textrm{and }c:=c_{21}^{21}-c_{11}^{22}\ ,
$$
we transform it into:
\\
$\hspace*{14.4em}
h_{12}\hspace*{0.9em} h_{11}\hspace*{0.9em} h_{22}
$
$$
\left(
\begin{array}{ccc}
0 & -b & c \\
0 & b & -c \\
-2c & 2a & 0 \\
-2b & 0 & 2a \\
\end{array} \right)
$$
It is clearly degenerate, and has rank 2 in general position.

This means we need to consider (\ref{eq:3.7}) in full generality:
for arbitrary $\Gamma_0$ and $\Gamma_1$.
((\ref{eq:3.10}) is (\ref{eq:3.7}) under assumption that $\Gamma_0=0$, which now has to be lifted.)\\
(\ref{eq:3.7}) involves $\widetilde{\mathcal{L}}_{V_2}\Gamma_0$, so we need to express $V_2$ from the
first equation of $\Gamma$-part of (\ref{sys:SYS}):
$$\mathcal{L}_{V_{1}}\Gamma_{0}+{\displaystyle\frac{\partial^{2}V_{2}}{\partial x^{2}}} =  0 \,.$$
Setting $\Gamma^k_{0ij}=:\gamma^k_{ij}$, it can be rewritten in index form as:
$$V_{2,ij}^l=\gamma_{ij}^{k}b^{l}_k-\gamma_{kj}^{l}b_i^{k}-\gamma_{ik}^{l}b_j^{k}=:v^{l}_{ij}$$
Actually, second derivatives of $V_2$ is all we need in (\ref{eq:3.7}), where they appear in
$$(\widetilde{\mathcal{L}}_{V_2}\Gamma_0)^l_{ij,p}
=-\gamma^k_{ij}V_{2,kp}^l+\gamma^l_{kj}V_{2,ip}^k+\gamma^l_{ik}V_{2,jp}^k\ ,$$
which we can now rewrite as:\\\\
$(\widetilde{\mathcal{L}}_{V_2}\Gamma_0)^l_{ij,p}=
-\gamma^k_{ij}(\gamma_{kp}^{s}b^{l}_s-\gamma_{sp}^{l}b_k^{s}-\gamma_{ks}^{l}b_p^{s})
+\gamma^l_{kj}(\gamma_{ip}^{s}b^{k}_s-\gamma_{sp}^{k}b_i^{s}-\gamma_{is}^{k}b_p^{s})$
\begin{equation}
\label{eq:LV2Gamma0}
\hspace*{14em}
+\gamma^l_{ik}(\gamma_{jp}^{s}b^{k}_s-\gamma_{sp}^{k}b_j^{s}-\gamma_{js}^{k}b_p^{s})
\end{equation}
One note about coefficients $\gamma$. Compatibility conditions (\ref{eq:compability})
in dimension $n=2$ become:
$$
\gamma_{1j}^{1}=-\gamma_{2j}^{2}\ .
$$
That leaves 4 independent coefficients: $\gamma_{11}^{2}$, $\gamma_{11}^{1}$, $\gamma_{12}^{1}$
and $\gamma_{22}^{1}$.

We consider (\ref{eq:3.7}) as $S(V_1)=0$ - linear operator acting on
$V_1$, and split the operator into two parts:
$S=S(\Gamma_0)+S(\Gamma_1)$. Matrix of $S(\Gamma_1)$ is calculated at the top of this section.\\
(\ref{eq:LV2Gamma0}) allows us to rewrite
$S(\Gamma_0) V_1=(\widetilde{\mathcal{L}}_{V_2}\Gamma_0)^l_{ij,p}-
(\widetilde{\mathcal{L}}_{V_2}\Gamma_0)^l_{pj,i}$
as:\\\\
$
( \gamma^k_{pj}\gamma_{ki}^{s}- \gamma^k_{ij}\gamma_{kp}^{s} ) b^{l}_s +
( \gamma^l_{pk}\gamma_{js}^{k}- \gamma^k_{pj}\gamma_{ks}^{l} ) b^{s}_i +$
\begin{equation}
\label{eq:comp:gamma0}
\hspace{9em}+( \gamma^l_{pk}\gamma_{si}^{k}- \gamma^l_{ik}\gamma_{sp}^{k} ) b^{s}_j+
( \gamma^k_{ij}\gamma_{ks}^{l}- \gamma^l_{ik}\gamma_{js}^{k} ) b^{s}_p
\end{equation}
Recall that dimension $n=2$, and indexes $1=i<p=2$ must therefore
stay fixed at $i=1$, $p=2$, while the remaining pair of indexes
take any values.  That turns (\ref{eq:comp:gamma0}) into a
system of 4 expressions indexed with $(j,l)$:\\\\
$(11)\hspace{2em}( \gamma^1_{22}\gamma_{11}^{2}-
\gamma^1_{12}\gamma_{12}^{2} ) b^{1}_1 +
             ( \gamma^2_{11}\gamma_{22}^{1}- \gamma^1_{12}\gamma_{12}^{2} ) b^{2}_2 + $\\\\
$\hspace*{17em}    +( \gamma^1_{2k}\gamma_{12}^{k}- \gamma^1_{1k}\gamma_{22}^{k} ) b^{2}_1 +
             ( \gamma^k_{21}\gamma_{k1}^{2}- \gamma^k_{11}\gamma_{k2}^{2} ) b^{1}_2$\\\\
$(22)\hspace{2em}( \gamma^2_{21}\gamma_{21}^{1}- \gamma^1_{22}\gamma_{11}^{2} ) b^{1}_1 +
             ( \gamma^1_{12}\gamma_{12}^{2}- \gamma^2_{11}\gamma_{22}^{1} ) b^{2}_2 + $\\\\
$\hspace*{17em}    +( \gamma^k_{22}\gamma_{k1}^{1}- \gamma^k_{12}\gamma_{k2}^{1} ) b^{2}_1 +
             ( \gamma^2_{2k}\gamma_{11}^{k}- \gamma^2_{1k}\gamma_{21}^{k} ) b^{1}_2\hfill$\\\\
$(12)\hspace{2em}2( \gamma^2_{2k}\gamma_{11}^{k}- \gamma^k_{21}\gamma_{k1}^{2} ) b^{1}_1 +
\hspace{5em}2( \gamma^2_{21}\gamma_{21}^{1}- \gamma^2_{11}\gamma_{22}^{1})b^{2}_1\hfill$\\\\
$(21)\hspace{8em}2( \gamma^1_{2k}\gamma_{21}^{k}- \gamma^1_{1k}\gamma_{22}^{k} ) b^{2}_2 +
\hfill 2( \gamma^1_{22}\gamma_{11}^{2}- \gamma^1_{12}\gamma_{12}^{2})b^{1}_2$\\\\
We use (\ref{comp:Gamma1}, $V_1$-hamiltonian) to go from
$b$-coefficients for $V_1$ to $h$-coefficients. Then the fact that
$V_2$ too is hamiltonian (second equation in $\omega$-part of
(\ref{sys:SYS})) follows automatically, as a short calculation
would show. Considered by itself, this system is degenerate.
Indeed, setting $$( \gamma^1_{22}\gamma_{11}^{2}-
\gamma^1_{12}\gamma_{12}^{2} )=:A, ( \gamma^1_{2k}\gamma_{12}^{k}-
\gamma^1_{1k}\gamma_{22}^{k} )=:B, ( \gamma^k_{21}\gamma_{k1}^{2}-
\gamma^k_{11}\gamma_{k2}^{2} )=:C\ ,$$
and using $h$-coefficients, the system's matrix becomes:
\\
$\hspace*{14em} h_{12}\hspace*{1.35em} h_{11}\hspace*{1.33em}
h_{22} $
$$
\left(
\begin{array}{ccc}

        0         &          -B         &          C         \\
        0         &           B         &         -C         \\
       -2C        &          2A         &          0         \\
       -2B        &           0         &         2A        \\
\end{array} \right)
$$\\
The determinant of this is identically zero.\\
Notice that the two matrices for $S(\Gamma_0)$ and
 $S(\Gamma_1)$ obtained so far look exactly the same, up to
capitalization of the entries' names. The matrix for the operator
$S=S(\Gamma_0)+S(\Gamma_1)$ is a sum of the two. Since it will
have the same structure as either of its degenerate summands, it is
also degenerate. It has rank 2 however, since it's lower right 2 X 2 block is:
$$\left(
\begin{array}{cc}
        2(a+A)      & 0        \\
        0       &     2(a+A) \\
\end{array} \right)$$
This is non-degenerate in general position, since:
$$
a+A=\gamma^1_{22}\gamma_{11}^{2}+\gamma^1_{12}\gamma_{11}^{1}-(c_{21}^{22}+c_{12}^{11})\neq0\ ,
$$
resulting in a 1-dimensional stabilizer at 1-jet.

Let us now consider the next, second jet of our connection. To calculate its stabilizer,
we need to solve the following equation from (\ref{sys:SYS}) for $V_4$ :
$${\mathcal{L}}_{V_{1}}\Gamma_{2}+\tilde{\mathcal{L}}_{V_{2}}\Gamma_{1}+
\tilde\mathcal{L}_{V_{3}}\Gamma_{0}+
{\displaystyle\frac{\partial^{2}V_{4}}{\partial x^{2}}}  =  0
$$
Its compatibility conditions are:
\begin{equation}
\label{eq:compV4}
({\mathcal{L}}_{V_{1}}\Gamma_{2}+\tilde{\mathcal{L}}_{V_{2}}\Gamma_{1}+
\tilde\mathcal{L}_{V_{3}}\Gamma_{0})_{ij,p}^l=  ({\mathcal{L}}_{V_{1}}\Gamma_{2}+\tilde{\mathcal{L}}_{V_{2}}\Gamma_{1}+
\tilde\mathcal{L}_{V_{3}}\Gamma_{0})_{pj,i}^l
\end{equation}
We will use the same strategy as in the previous section
to prove that in this case stabilizer is trivial. Namely we will obtain a connection
2-jet, for which the above equation will be a non-degenerate homogeneous linear system.
We set $\Gamma_0=\Gamma_1=0$. This implies $V_2=V_3=0$, hence hamiltonian, so that
$\omega$-part of (\ref{sys:SYS}) is true for any hamiltonian $V_1$.
That simplifies (\ref{eq:compV4}) to:
\begin{equation}
\label{eq:compV4:restricted}
({\mathcal{L}}_{V_{1}}\Gamma_{2})_{ij,p}^l
=({\mathcal{L}}_{V_{1}}\Gamma_2)_{pj,i}^l
\end{equation}
We introduce notation for coefficients of $\Gamma_2$ :
$${\Gamma_2}_{ij}^l=\sum_{s,t=1}^{2}d^l_{ijst}x^s x^t\ ,\ d^l_{ijst}=d^l_{jits}$$
Compatibility with $\omega$ (\ref{eq:compability}) impose these restrictions on $d$ in dimension 2:
$$
d^2_{2\alpha ij}=-d^1_{1\alpha ij}
$$
There are thus 4 families of independent coefficients: $d^2_{11ij}$, $d^1_{11ij}$, $d^1_{12ij}$ and $d^1_{22ij}$.
With these, $$({\mathcal{L}}_{V_{1}}\Gamma_{2})_{ij,p}^l=
2d^l_{ijkt} b^k_p x^t+2d^l_{ijkp} b^k_t x^t-2d^k_{ijpt} b^l_k x^t+2d^l_{kjpt} b^k_i x^t+2d^l_{ikpt} b^k_j x^t$$
( $b_k^l$ are still coefficients of $V_1$, as in section (\ref{sec:Proof}), and (\ref{eq:compV4:restricted})
( with $i=1,p=2$) is:\\\\
$(d^l_{kj2t}-d^l_{2jkt})b^k_1 + (d^l_{1jkt}-d^l_{kj1t})b^k_2 + (d^l_{1jk2}-d^l_{2jk1})b^k_t +$
\\\\
$\hspace*{17em}
(d^k_{2j1t}-d^k_{1j2t})b^l_k + (d^l_{1k2t}-d^l_{2k1t})b^k_j = 0 $\\\\
With the triple of indexes $(j,l,t)$ arbitrary, we have system of 8 equations in 4 variables: the coefficients of $V_1$.
This is the system, equations are labelled by this index triple:\\\\
$
(111)
\hspace{2em}2(d^1_{1112}-d^1_{1211})b_1^1+(d^1_{1112}-d^1_{1211})b_2^2$\\
$\hspace*{17em}+(d^1_{1112}-d^1_{1211})b_1^2+  (d^2_{1211}-d^2_{1112})b_2^1=0$\\\\
$
(221)
\hspace{2em}2(d^2_{1212}-d^2_{2211})b_1^1+(d^2_{1212}-d^2_{2211})b_2^2+$\\
$\hspace*{12em}(d^2_{1222}-d^2_{2221}+d^1_{2211}-d^1_{1212})b_1^2+(d^2_{1112}-d^2_{1211})b_2^1=0
$\\\\
$
(121)
\hspace{2em}3(d^2_{1112}-d^2_{1211})b_1^1+\hspace{3em}(d^2_{1122}-d^2_{2211}+d^1_{2111}-d^1_{1112})b_1^2
\hspace{3em}=0
$\\\\
$
(211)
\hspace{2em}(d^1_{1212}-d^1_{2211})b_1^1+2(d^1_{1212}-d^1_{2211})b_2^2+$\\
$\hspace*{12em}(d^1_{1222}-d^1_{2221})b_1^2+ (d^1_{1112}-d^1_{2111}+d^2_{2211}-d^2_{1221})b_2^1=0
$\\\\
$
(112)
\hspace{2em}(d^1_{1122}-d^1_{1212})b_1^1+2(d^1_{1122}-d^1_{1212})b_2^2+$\\
$\hspace*{12em}(d^1_{1222}-d^1_{2221})b_1^2+(d^1_{1112}-d^1_{2111}+d^2_{1212}-d^2_{1122})b_2^1=0\\
$\\\\
$
(222)
\hspace{2em}(d^2_{1222}-d^2_{2221})b_1^1+2(d^2_{1222}-d^2_{2221})b_2^2+$\\
$\hspace*{17em}(d^1_{2221}-d^1_{1222})b_1^2+(d^2_{1122}-d^2_{2211})b_2^1=0\\
$\\\\
$
(122)
\hspace{2em}2(d^2_{1122}-d^2_{1212})b_1^1+(d^2_{1122}-d^2_{1212})b_2^2+$\\
$\hspace*{12em}(d^2_{1222}-d^2_{2221}+d^1_{1212}-d^1_{1122})b_1^2+(d^2_{1112}-d^2_{2111})b_2^1=0\\
$\\\\
$
(212)
\hspace{5.5em}3(d^1_{1222}-d^1_{2212})b_2^2+\hspace{4em}(d^2_{2221}-d^2_{1222}+d^1_{1122}-d^1_{2211})b_2^1=0$
\\\\
Setting:
$$
a=(d^1_{1112}-d^1_{1211}), e=(d^1_{1212}-d^1_{1122}), g=(d^1_{1222}-d^1_{2221}), h=(d^1_{1122}-d^1_{1212})\ ,$$
$$b=(d^2_{1211}-d^2_{1112}), c=(d^2_{1212}-d^2_{2211}), d=(d^2_{1222}-d^2_{2221}), f=(d^2_{1122}-d^2_{1212})
\ ,$$
we see the system take form:
\\
$\hspace*{2.7em} b_1^1\hspace*{1.8em} b_2^2\hspace*{3em} b_1^2\hspace*{4.1em} b_2^1
\hspace*{5.8em} h_{12}\hspace*{2.7em} h_{11}\hspace*{3.3em}
h_{22} $
$$\left(
\begin{array}{cccc}
2a & a & a   & b \\
2c & c & d+e &-b \\
-3b& 0 &f+c-a& 0 \\
-e & -2e & g & a-c\\
h  & 2h &  g & a-f\\
d  & 2d & -g & f+c\\
2f &  f & d-h& -b \\
 0 & 3g &  0 &h-e-d\\
\end{array} \right)
=
\left(
\begin{array}{ccc}
a & -a   & b \\
c & -d-e &-b \\
-3b&   a-f-c&  0 \\
e & -g & a-c\\
-h  &  -g & a-f\\
-d & g & f+c\\
f & h-d& -b \\
-3g & 0 &h-e-d\\
\end{array} \right)
$$
This is non-degenerate for a generic connection.
For example, if $d^2_{1112}=d^2_{1122}=1$, the rest is null,
then $f=1$, $b=-1$, all others zero, and the system is:
$$\left(
\begin{array}{ccc}
 & & -1 \\
& &   1 \\
3&  -1 &   \\
 &  & \\
 &  & -1  \\
  &  & 1 \\
1 &  & 1 \\
  &  &   \\
\end{array} \right)$$

Now we can summarize what we know about exceptional stabilizers:
\begin{prop}\label{except:stab:dim}
The stabilizer of a k-jet of a generic connections\\
for $n=2$ is equal to $G_1/G_2$ for $k=0$, is 1-dimensional for $k=1$ ,\\
and is trivial for $k\geq 2$ .
\end{prop}
\section{Poincar\'{e} series}
\label{sec:Proof} Here we will calculate the Poincar\'e series of
$\mathcal{M}$, the moduli space of Fedosov structures:
$$
p_{\Phi}(t) = \dim \mathcal{M}_{0}+\sum_{k=1}^{\infty}
(\dim\mathcal{M}_{k}- \dim \mathcal{M}_{k-1})t^{k}
$$
To obtain $\dim\mathcal{M}_{k}$, we need to discuss $\mathcal{F}_{k}$ first.
In particular, we need to know how many different local symplectic structures are there.
More precisely, we want to find the dimension of the space of $k$-jets of non-degenerate
closed 2-forms at a point. Non-degeneracy is an open condition and does not affect dimension.
Closedness locally is equivalent to exactness. For a symplectic form $\omega$:
$$\omega=d \alpha\ ,$$
for some 1-form $\alpha$ defined up to $\nabla f$, a gradient of a function, that function in its
turn is defined up to a constant.
We have the following exact sequence:
$$
0 \longrightarrow \mathbb{R} \longrightarrow C^\infty(\mathbb{R}^{2n})
\stackrel{d^0}{\longrightarrow} \Lambda^1(\mathbb{R}^{2n})
\stackrel{d^1}{\longrightarrow} d\Lambda^1(\mathbb{R}^{2n}) \longrightarrow 0\ ,
$$
which descends to jets:
$$
0 \longrightarrow \mathbb{R} \longrightarrow j^{l+2}(C^\infty(\mathbb{R}^{2n}))
\stackrel{d^0}{\longrightarrow} j^{l+1}(\Lambda^1(\mathbb{R}^{2n})) \stackrel{d^1}{\longrightarrow}
j^{l}(d\Lambda^1(\mathbb{R}^{2n}))
\stackrel{d^2}{\longrightarrow} 0
$$
It follows that:
$$
\dim[j^{l}d\Lambda^1(\mathbb{R}^{2n})]=
\dim[j^{l+1}\Lambda^1(\mathbb{R}^{2n})]-\dim[j^{l+2}(C^\infty(\mathbb{R}^{2n}))]+\dim\mathbb{R}
$$
We are interested in 0-jets since higher jets of $\omega$ are determined by
the connection part $\Gamma$ of a given Fedosov structure $\Phi$ through
compatibility condition (\ref{eq:compability}), see Theorem 4.5 (2) p.124 in \cite{GRS}.
$$
\dim[j^{0}d\Lambda^1(\mathbb{R}^{2n})]=
\dim[j^{1}\Lambda^1(\mathbb{R}^{2n})]-\dim[j^{2}(C^\infty(\mathbb{R}^{2n}))]+1=$$
$$
=2n{2n+1 \choose 2n}-{2n+2 \choose 2n}+1=\frac{2n(2n-1)}{2}
$$
Each $\omega$ is compatible with (or preserved by) all $\Gamma$,
such that $\omega_{i\alpha}\Gamma^\alpha_{jk}$
is completely symmetric in $i$, $j$, $k$, cf. the last paragraph on p.110 in \cite{GRS}.
At 0-jet of Fedosov structure $\Phi_0=(\omega_0,\Gamma_0)$ there are ${2n+3-1 \choose 2n-1}$ of those,
hence:
$$
\dim\mathcal{F}_0=\dim\{\textrm{all }\omega_0\}\dim\{\textrm{all compatible }\Gamma_0\}=
\frac{2n(2n-1)}{2}{2n+2 \choose 2n-1}
$$
For other $\mathcal{F}_k$'s we must remember that each $\Gamma_{jk}^{i}$ is a homogeneous polynomial
of degree $k$ in $2n$ variables:
$$
\dim\mathcal{F}_k=\frac{2n(2n-1)}{2}{2n+2 \choose 2n-1}\sum_{m=0}^k{2n+m-1 \choose 2n-1}=
\frac{2n(2n-1)}{2}{2n+2 \choose 2n-1}{2n+k \choose 2n}
$$
Next, we need to know orbit dimensions.
$G_1/G_3$ acts on $\Phi_0$ non-trivially, i.e. both first and second component
of generating vector field $V_1$ and $V_2$ are acting. The stabilizer $G_{\Phi_0}$
is determined by an arbitrary hamiltonian $V_1$:
$$
\dim\mathcal{O}_0=\dim\{(V_1,V_2)\}-\dim \texttt{sp}(2n)
$$
$$
=2n\sum_{m=1}^2{2n+m-1 \choose 2n-1}-\frac{2n(2n+1)}{2}=n((2n+1)^2-2)=n(4n^2+4n-1)
$$
$V_1$,..., $V_{k+2}$ act on $\Phi_k$:
$$
\dim\mathcal{O}_k=\dim\{(V_1,\ldots,V_{k+2})\}-1\cdot\delta_{2n}^2\delta_k^1
$$
(Kronecker symbol $\delta$ is needed here to take care of exceptional dimension two.)
$$
=2n\sum_{m=1}^{k+2}{2n+m-1 \choose 2n-1}-\delta_{2n}^2\delta_k^1=2n\left[{2n+k+2 \choose 2n}-1\right]
-\delta_{2n}^2\delta_k^1
$$
This gives us dimension of moduli space of $k$-jets:
$$
\dim \mathcal{M}_0=\frac{2n(2n-1)}{2}{2n+2 \choose 2n-1}-n((2n+1)^2-2)=
\frac{n[8n(2n^2-1)(n+1)+11]}{6}
$$
$$
\dim \mathcal{M}_k=\dim\mathcal{F}_k-\dim\mathcal{O}_k
$$
$$
=\frac{2n(2n-1)}{2}{2n+2 \choose 2n-1}\sum_{m=0}^k{2n+m-1 \choose 2n-1}-
2n\sum_{m=1}^{k+2}{2n+m-1 \choose 2n-1}+\delta_{2n}^2\delta_k^1
$$
$$
=\frac{2n(2n-1)}{2}{2n+2 \choose 2n-1}{2n+k \choose 2n}-2n\left[{2n+k+2 \choose 2n}-1\right]
+\delta_{2n}^2\delta_k^1\ ,k\geq 1
$$
We will have to write constant and linear terms of Poincar\'e series separately because
they contain $\mathcal{M}_0$. The linear coefficient is:
$$
\dim \mathcal{M}_1-\dim \mathcal{M}_0=\frac{n(2n+1)}{3}\left[4n^4+2n^3-6n^2-4n-3\right]+\delta_{2n}^2
$$
The common term will have this coefficient:
$$
\dim \mathcal{M}_k-\dim \mathcal{M}_{k-1}=\frac{2n(2n-1)}{2}{2n+2\choose 2n-1}{2n+k-1\choose 2n-1}-
2n{2n+(k+2)-1\choose 2n-1}
$$
$$
=4n{2n+2\choose 2n-2}{2n+k-1\choose 2n-1}-2n{2n+k+1\choose 2n-1}-\delta_{2n}^2\delta_{k}^2\ ,k\geq2
$$
We have:
$$
p_{\Phi}(t) = \dim \mathcal{M}_{0}+\sum_{k=1}^{\infty}
(\dim\mathcal{M}_{k}- \dim \mathcal{M}_{k-1})t^{k}
$$
$$
=\frac{n[8n(2n^2-1)(n+1)+11]}{6}+\frac{n(2n+1)}{3}\left[4n^4+2n^3-6n^2-4n-3\right]t+(t-t^2)\delta_{2n}^2+
$$
$$
{\displaystyle
+2n\sum_{k=2}^{\infty}\;\biggl[2{2n+2 \choose
4} {2n+k-1\choose2n-1} - {2n+k+1 \choose 2n-1}\biggr]t^{k}
}
$$
\begin{prop}
The Poncar\'e series $p_\Phi(t)$ is a rational function.
Namely,
$$
p_\Phi(t)=\frac{n(20n^2+8n+11)}{6}-\frac{n(2n+1)}{3}\left[4n^4+2n^3+2n^2-4n+3\right]t$$
$$
+(t-t^2)\delta_{2n}^2+2nD_\Phi\left(\frac{1}{1-t}\right)
$$
where $D_\Phi$ is a differential operator of order $2n-1$ :
$$D_{\Phi}=2{2n+2 \choose 4}
{2n+t\frac{d}{dt}-1\choose2n-1} - {2n+t\frac{d}{dt}+1 \choose 2n-1}$$
with $${2n+t\frac{d}{dt}-1 \choose 2n-1}=\frac{1}{(2n-1)!}(t\frac{d}{dt}+1)\ldots(t\frac{d}{dt}+2n-1)\ ,$$
$${2n+t\frac{d}{dt}+1 \choose 2n-1}=\frac{1}{(2n-1)!}(t\frac{d}{dt}+3)\ldots(t\frac{d}{dt}+2n+1) \ .$$
\end{prop}
{\bf Proof}\quad Indeed, denote
$$\varphi_{m}(t)=\sum_{k=0}^{\infty}k^{m}t^{k}\ ,\qquad m\in\mathbb{Z}_{+}\ ,$$
then
$$\varphi_{m}(t)=\sum_{k=0}^{\infty}k^{m-1}kt^{k-1}t=t\left(\sum_{k=0}^{\infty}k^{m-1}t^{k}\right)'
=\left(t\frac{d}{dt}\right)\varphi_{m-1}(t)\quad \mathrm{for}\ m\in\mathbb{N}\,.$$
Thus
$$\varphi_{m}(t)=\left(t\frac{d}{dt}\right)^{m}\varphi_{0}(t)=
\left(t\frac{d}{dt}\right)^{m}\left(\frac{1}{1-t}\right)\,.$$
Hence,
$$\sum_{k=0}^{\infty}\;\biggl[2{2n+2 \choose 4} {2n+k-1\choose2n-1} - {2n+k+1 \choose 2n-1}\biggr]t^{k}$$
$$=\left[2{2n+2 \choose 4}{2n+t\frac{d}{dt}-1\choose2n-1} - {2n+t\frac{d}{dt}+1 \choose 2n-1}\right]
\left(\frac{1}{1-t}\right)\,.$$
We have:
$$
\sum_{k=2}^{\infty}\;\biggl[2{2n+2 \choose 4} {2n+k-1\choose2n-1} - {2n+k+1 \choose 2n-1}\biggr]t^{k}
$$
$$=\sum_{k=0}^{\infty}\;\biggl[2{2n+2 \choose 4} {2n+k-1\choose2n-1} - {2n+k+1 \choose 2n-1}\biggr]t^{k}$$
$$
-{2n+2 \choose 3}(2n^2-n-1)t-{2n+1 \choose 2}\frac{2n^2+n-4}{3}
$$
So
$$p_{\Phi}(t)=\frac{n[8n(2n^2-1)(n+1)+11]}{6}+\frac{n(2n+1)}{3}\left[4n^4+2n^3-6n^2-4n-3\right]t
$$
$$
+(t-t^2)\delta_{2n}^2+2n\Bigg\{\sum_{k=0}^{\infty}\;\biggl[2{2n+2 \choose 4} {2n+k-1\choose2n-1} - {2n+k+1 \choose 2n-1}\biggr]t^{k}
$$
$$
-{2n+2 \choose 3}(2n^2-n-1)t-{2n+1 \choose 2}\frac{2n^2+n-4}{3}\Bigg\}
$$
$$
=\frac{n(20n^2+8n+11)}{6}-\frac{n(2n+1)}{3}\left[4n^4+2n^3+2n^2-4n+3\right]t
$$
$$
+(t-t^2)\delta_{2n}^2+2nD_\Phi\left(\frac{1}{1-t}\right)$$
\hfill$\Box$


\begin{thebibliography}{99}

\bibitem[A]{Arn} Arnold, V.I., Mathematical Problems in Classical Physics in:\\
\textit{Trends and Perspectives in Applied Mathematics,\\
Applied Mathematics Series} vol.100, Editors: F.John, J.E.Marsden, L.Sirovich;
N.Y. Springer 1999, pp.1-20.

\bibitem[D]{Dubr} Dubrovskiy, S. Moduli space of symmetric connections,
\textit{``Representation Theory, Dynamical Systems, Combinatorial and Algorithmic Methods.
Part 7''} (A.M.Vershik ed.) \textit{Zapiski Nauchnyh Seminarov POMI} {\bf 292}(2002), 22-39.
Available electronically: \texttt{http://www.pdmi.ras.ru/znsl/2002/v292.html}


\bibitem[G]{Gersh} Gershkovich, V. Ya. On normal form of distribution jets. \textit{Topology and geometry-Rohlin Seminar, Lecture Notes in Math.}, 1346, pp.77-98, Springer, Berlin, 1988.

\bibitem[VG]{Versh:Gersh} A.Vershik, V.Gershkovich,
Estimation of the functional dimension of the orbit space of germs of distributions in general position. (Russian)
\textit{Mat. Zametki} 44 (1988), no. 5, 596-603, 700; translation in \textit{Math. Notes} 44 (1988), no. 5-6, 806-810 (1989)

\bibitem[W]{Weym} P.Magyar, J.Weyman, A.Zelevinsky, Multiple flag varieties of finite type.
\textit{Adv. Math.} 141 (1999), no. 1, 97--118.

\bibitem[Sh]{Shmel:struc} Shmelev, A.S., On Differential Invariants of Some Differential-Geometric Structures,
\textit{Proceedings of the Steklov Institute of Mathematics}, 1995, vol.209, pp.203-234.

\bibitem[Sh2]{Shmel:mod} Shmelev, A.S., Functional moduli of germs of Riemannian metrics (Russian), \textit{Funktsional. Anal. i Prilozhen.}, 31 (1997), no. 2, pp.58--66, 96; translation in \textit{Funct. Anal. Appl.} 31 (1997), no. 2, pp.119--125

\bibitem[T]{T} Tresse, A., Sur les Invariants Diff\'erentiels des Groupes Continus des\\
Transformations, \textit{Acta Mathematica}, 1894, vol.18, pp.1-88.

\bibitem[GRS]{GRS} Gelfand I.M., Retakh, V., Shubin, M.A, Fedosov Manifolds,\\



\end{thebibliography}
\end{document}